\documentstyle[a4,amsfonts,12pt]{article}

\newcommand{\mb}{\mbox}
\newcommand{\beq}{\begin{equation}}
\newcommand{\eeq}{\end{equation}}
\newcommand{\ueberschrift}{\bigskip\goodbreak\noindent\bigskip}
\newcounter{theabsatz}
\newcommand{\absatz}[1]{\stepcounter{theabsatz} \ueberschrift
               {\large \bf \arabic{theabsatz}. {#1}} \setcounter{equation}{0}}

 \newtheorem{theor}{Theorem}
 
 \newtheorem{lem}{Lemma}

\begin{document}
\newcommand{\ext}{\mbox{ext\,}}
\newcommand{\diam}{\mbox{diam\,}}

\parindent 10 pt
\parskip 8pt plus 4pt
\jot 10pt

\abovedisplayskip 8pt plus 1pt \belowdisplayskip 8pt plus 1pt

\newcommand{\C}{{\mathbb{C}}}
\newcommand{\CC}{{\overline{\mathbb{C}}}}
\newcommand{\DD}{{\overline{\mathbb{D}}}}
\newcommand{\D}{{\mathbb{D}}}
\newcommand{\R}{{\mathbb{R}}}
\newcommand{\He}{{\mathbb{H}}}
\newcommand{\T}{{\mathbb{T}}}
\newcommand{\N}{{\mathbb{N}}}
\newcommand{\PP}{{\mathbb{P}}}
\newcommand{\p}{\preceq}
\newcommand{\s}{\succeq}
\newcommand{\En}{{\mathcal E}_n}
\newcommand{\ov}{\overline}
\newcommand{\de}{\delta}
\newcommand{\Ga}{\Gamma}
\newcommand{\ga}{\gamma}
\newcommand{\la}{\lambda}
\newcommand{\kap}{\mb{cap}}
\newcommand{\be}{\beta}
\newcommand{\Om}{\Omega}
\newcommand{\om}{\omega}
\newcommand{\al}{\alpha}
\newcommand{\ve}{\varepsilon}

\begin{center}
{\large \bf    Polynomial approximation 
on a compact subset of the real line}
\\[3ex] { Vladimir Andrievskii}\\[3ex]

%{\it  Department of Mathematical Sciences, Kent State University,
% Kent, OH 44242}\\[4ex]
{\bf Abstract}
\end{center}
\vspace*{0.5cm}
 We prove an analogue of the classical Bernstein theorem concerning the
rate of polynomial approximation of piecewise analytic functions
on a compact subset of the real line.\\



\absatz{Introduction and the main result}

Let $E\subset \R$ be a compact subset of the real line $\R$ and
let $\PP_n$ be the set of all (real) polynomials of degree at most
$n\in\N:=\{1,2,\ldots\}$. Also, let for $x_0\in E$ and $\al>0$,
 $$\En(|x-x_0|^\al,E):=\inf_{p\in\PP_n}\sup_{x\in E}||x-x_0|^\al-p(x)|.
 $$
  The starting point of our analysis is the
	classical
 Bernstein theory \cite{ber1, ber2, ber3}. According to this theory,
for any $x_0\in(-1,1)$ and $\al>0$, where $\al $ is not an even
integer, there exists a  finite  nonzero limit
 $$\sigma_\al:=\lim_{n\to\infty}n^\al\En(|x-x_0|^\al,[-1,1]).
 $$
 The question as to what happens to the best polynomial approximations
 for a general set $E\subset\R$ 
  is investigated in  monographs \cite{vas} and \cite[Chapter 10]{tot06} where the
 reader can also find  a comprehensive survey of this subject.

Now, we consider $E$ to be a set in the complex plane $\C$ and use
 the notions  of potential theory in the plane (see \cite{ran,
 saftot} for details).
  Let $E$ be non-polar, i.e., be of positive (logarithmic) capacity
 cap$(E)>0$ and let $g_{\CC\setminus E}(z)=g_{\CC\setminus E}(z,\infty),
 z\in\CC\setminus E$
 be the Green
 function of $\CC\setminus E$ with pole at infinity,
  where $\CC:=\C\cup\{\infty\}$ is the extended complex plane.

Our main objective is to prove the following result.
\begin{theor}\label{th1}
Let $x_0\in E$. If for some $\al>0$,
\beq\label{1.1}
\limsup_{n\to\infty}n^\al\En(|x-x_0|^\al,E)>0,
\eeq
then
\beq\label{1.2}\sup_{z\in\C\setminus E}\frac{g_{\CC\setminus E}(z)}{|z-x_0|}<
\infty.
\eeq
\end{theor}
Comparing Theorem \ref{th1} with \cite[Corollary 1]{and09} we obtain the following
result.
\begin{theor}
 Let $x_0\in E$. Then for any $\al>0$, which is not even integer,
$$
\liminf_{n\to\infty}n^\al\En(|x-x_0|^\al,E)>0
$$
if and only if (\ref{1.2}) holds.
\end{theor}
For the geometry of $E$ satisfying (\ref{1.2}), we
 refer the reader to \cite{cartot, tot06, cargar, and08} and the
 many references therein. 

\absatz{Auxiliary results}

In this section we assume that
$$
E=\bigcup_{j=1}^m[a_j,b_j],\quad x_0\in \bigcup_{j=1}^m (a_j,b_j)=:
\mb{ Int}(E),
$$
where $-1=a_1<b_1<a_2<\ldots <b_{m-1}<a_m<b_m=1, m>1.$

It is known (for example, see  \cite[pp. 224-226]{wid}, \cite[pp.
409--412]{andbla} or \cite{and01}) that 
there exists a conformal mapping $w=F(z)=F_E(z)$ of the upper
half-plane $\He:=\{z:\, \Im z>0\}$ onto   the
domain 
$$G=G_E=\{ z:\, 0<\Re z<\pi,\Im z>0\}\setminus
\cup_{j=1}^{m-1}[u_j,u_j+iv_j],$$ 
where
$0=:u_0<u_1<u_2<\ldots<u_{m-1}<u_m:=\pi$ and $v_j>0,\,
j=1,\ldots,m-1$, which can be extended continuously to
 $\ov{\He}$ satisfying the following
boundary correspondence 
$$F(\infty)=\infty,\,
F((-\infty,-1])=\{z:\, \Re z=0,\Im z\ge 0\},$$ $$
F([1,\infty))=\{z:\, \Re z=\pi,\Im z\ge 0\},$$
$$F([a_j,b_j])=[u_{j-1},u_j],\quad j=1,\ldots,m,$$
$$F([b_j,a_{j+1}])=[u_j,u_j+iv_j],\quad j=1,\ldots,m-1.$$ 
Moreover,
$$
g_{\CC\setminus E}(z)=\Im(F(z)),\quad z\in \ov{\He},
$$ 
$$
\pi\mu_E([a,b])=|F([a,b]\cap E)|,\quad [a,b]\subset [-1,1],$$
where $\mu_E$ is the equilibrium measure for $E$ and $|S|$ means
the linear Lebesgue measure (length) of $S\subset\C$.

Since, by the reflection principle, $F$ can be extended analytically
in a neighborhood of $x_0$, we can consider function
$$
h(x_0)=h_E(x_0):=F'(x_0)=\pi\om_E(x_0),
$$
where $\om_E(x_0)$ is the density of $\mu_E$ at $x_0$.

Our interest in this function lies in the fact that by the Vasiliev - Totik theorem
(see  \cite[p. 163, Corollary 11]{vas} and \cite[Theorem 10.5]{tot06}) for $\al>0$,
\beq\label{2.3}
\lim_{n\to\infty}n^\al\En(|x-x_0|^\al, E)=h(x_0)^{-\al}\sigma_\al.
\eeq
We set a sufficiently large constant $c>2$ such that
\beq\label{2.2}
h(x_0)\le c\quad \mb{and}\quad \kap(E)>\frac{1}{c}\, .
\eeq
In what follows,  we denote by $c_1,c_2,\ldots$ positive
constants which depend only on $c$.

By virtue of \cite[p. 39, (1.8)]{and01}, for $j=1,\ldots,m-1,$
\begin{eqnarray*}
v_j&\le&\log\left(\frac{1+(1-2\kap(E))^{1/2}}{1-(1-2\kap(E))^{1/2}}\right)
\\
&\le&
\log\left(\frac{\sqrt{c}+\sqrt{c-2}}{\sqrt{c}-\sqrt{c-2}}\right)
<\log(2c).
\end{eqnarray*}
For the fixed point $z_0:=x_0+2ce^{4\pi}i\in\He$, $w_0:=F(z_0)$,
and $\D:=\{z:|z|<1\}$, we have
\begin{eqnarray}
\Im(w_0)&=& g_{\CC\setminus E}(z_0)\ge g_{\CC\setminus \ov{\D}}(z_0)
=\log|z_0|\ge\log|z_0-x_0|
\nonumber\\
\label{2.7}
&=& 4\pi+\log(2c)=:c_1>v_j+4\pi,\quad j=1,\ldots,m-1.
\end{eqnarray}
Meanwhile, by virtue of \cite[p. 53]{saftot} and (\ref{2.2}),
\begin{eqnarray}
\Im(w_0)&=& g_{\CC\setminus E}(z_0)
=\int\log|z_0-x|d\mu_E(x)-\log\kap(E)
\nonumber\\
\label{2.17}
&\le&\log(|z_0-x_0|+2)+\log c<\log(4c^2)+4\pi<2c_1.
\end{eqnarray}
We use the notion of the {\it module
of a family of curves}.
We refer to
  \cite[Chapter 4]{ahl}, \cite[Chapter 9]{pom} or
\cite[pp. 341-360]{andbla} for the definition and basic properties of
the module (such as conformal invariance, comparison principle,
composition laws, etc.) We use this properties without further 
citation.
The curves under consideration are
crosscuts (see \cite{pom}) either of $\He$ or of $G$.

For $z_1,z_2\in\He$ denote by $\Ga(x_0,z_1;z_2,\infty;\He)$
the family of all crosscuts of $\He$ that separate  points $x_0$ and
$z_1$ from $z_2$ and $\infty$ in $\He$. Note that if $|z_1-x_0|<|z_2-x_0|,$ then
  for the module of $\Ga:=\Ga(x_0,z_1;z_2,\infty;\He)$ we have
\beq\label{2.8}
0\le m(\Ga)- \frac{1}{\pi}\log\left|
\frac{z_2-x_0}{z_1-x_0}\right|\le 2.
\eeq
Indeed, let $\Ga_0\subset\Ga$ denote the family of all crosscuts of
$A:=\{z\in\He:|z_1-x_0|<|z-x_0|<|z_2-x_0|\}$ separating  the 
circular parts of the  boundary $\partial A$ in $A$. 
Then, by \cite[p. 347, Example 1.8]{andbla},
$$
m(\Ga)\ge m(\Ga_0)=\frac{1}{\pi}\log\left|\frac{z_2-x_0}{z_1-x_0}\right|,
$$
which implies the left hand side of (\ref{2.8}).

Next, consider the metric
$$
\rho(z):=\left\{\begin{array}{ll} (\pi|z-x_0|)^{-1}&\mb{ if }
z\in\He, e^{-\pi}|z_1-x_0|\le|z-x_0|\le e^\pi|z_2-x_0| ,\\[2ex] 0&\mb{
elsewhere in $\C$}.\end{array}\right.
$$
The analysis analogous to the proof of
\cite[p. 349, (1.11)]{andbla} yields
$$
\int_\ga\rho(z)|dz|\ge 1,\quad \ga\in\Ga,
$$
i.e., $\rho$ is {\it admissible} (in the L-definition) of $m(\Ga)$.
Therefore,
$$
m(\Ga)\le\int\rho(z)^2dm(z)=\frac{1}{\pi}\log\left|
\frac{z_2-x_0}{z_1-x_0}\right| +2,
$$
where $dm$ means integration with respect to the 
two-dimensional Lebesgue measure
(area),
which implies the right-hand side inequality in (\ref{2.8}).

In a similar way, for $\eta_0:= F(x_0)$ and $w_1,w_2\in G$, we introduce the family
$\Ga(\eta_0,w_1;w_2,\infty;G)$
of all crosscuts of $G$ that separate  points $\eta_0$ and
$w_1$ from $w_2$ and $\infty$ in $G$.

For $r>0$, denote by $\ga(r)=\ga(\eta_0,G,r)\subset\{w\in G:|w-\eta_0|=r\}$
the crosscut of $G$ which has nonempty intersection with the ray
$\{w\in \He: \Re(w)=\eta_0\}$. For $0<r<R$, denote by $Q(r,R)=Q(\eta_0,G,r,R)\subset
G$ the bounded simply connected domain whose boundary consists of $\ga(r),\ga(R)$,
and two connected parts of $\partial G$. Let $m(r,R)=m(\eta_0,G,r,R)$ be the module of the family 
$\Ga(r,R)=\Ga(\eta_0,G,r,R)$ of all
crosscuts of $Q(r,R)$ which separate  circular arcs $\ga(r)$ and $\ga(R)$ in $Q(r,R)$.
\begin{lem}\label{lem2.1}
For $0<r<R\le R_0:=|w_0-\eta_0|$,
\beq\label{2.9}
m(r,R)\le\frac{1}{\pi}\log\frac{R}{r}+c_2,\quad c_2:=\frac{1}{\pi}\log\frac{c^2}{c_1}+7.
\eeq
\end{lem}
{\bf Proof}. From the right-hand side of (\ref{2.8}) we conclude 
that for $\Ga_1:=\Ga(x_0,z;z_0,\infty;\He)$ and $z\in\He$ with
$|z-x_0|<|z_0-x_0|$, 
$$m(\Ga_1)\le
\frac{1}{\pi}\log\left|\frac{z_0-x_0}{z-x_0}\right|+2.
$$
While, for $\Ga_1':=F(\Ga_1)=\Ga(\eta_0,w;w_0,\infty;G)$ and  $w:=F(z)$,
 where $|w-\eta_0|$ is sufficiently small, we have
\begin{eqnarray*}
m(\Ga_1')&\ge& m(|w-\eta_0|,R_0)\ge m(|w-\eta_0|,r) + m(r,R) +m(R,R_0)\\
&\ge& m(\{\{\xi\in\He:|\xi-x_0|=t\}:|w-\eta_0|<t<r\})
+m(r,R)\\
&&+m(\{\{\xi\in\He:|\xi-x_0|=t\}:R<t<R_0\})\\
 &=& \frac{1}{\pi}\log\frac{r}{|w-\eta_0|}+ m(r,R)
+\frac{1}{\pi}\log\frac{R_0}{R}\, .
\end{eqnarray*}
Since $m(\Ga_1)=m(\Ga_1')$, comparing above inequalities, we obtain
\begin{eqnarray*}
m(r,R)&\le&
\frac{1}{\pi}\log\left|\frac{z_0-x_0}{z-x_0}\right|+2
-\frac{1}{\pi}\log\frac{r}{|w-\eta_0|}-\frac{1}{\pi}\log\frac{R_0}{R}\\
&=& 
\frac{1}{\pi}\log\left|\frac{w-\eta_0}{z-x_0}\right|
+\frac{1}{\pi}\log\frac{|z_0-x_0|}{R_0} +2
+\frac{1}{\pi}\log\frac{R}{r}\, .
\end{eqnarray*}
Taking limit as $z\to x_0$ we have
$$
m(r,R)\le\frac{1}{\pi}\log |F'(x_0)|
+\frac{1}{\pi}\log\frac{|z_0-x_0|}{R_0} +2
+\frac{1}{\pi}\log\frac{R}{r}
$$
which,  together with (\ref{2.2}) and (\ref{2.7}), yields  (\ref{2.9}).

\hfill$\Box$

\begin{lem}\label{lem2.2}
For $j=1,\ldots, m-1$, we have
\beq\label{2.10}
\frac{v_j}{|\eta_0-u_j|}< 
\frac{\sqrt{(\eta_0-u_j)^2+v_j^2}}{|\eta_0-u_j|}\le
c_3:= 4\exp(3\pi c_2).
\eeq
\end{lem}
{\bf Proof}. 
Let  $r_j:=|\eta_0-u_{j}|, R_j:=\sqrt{|\eta_0-u_{j}|^2+v_{j}^2}.$
By Lemma \ref{lem2.1}
\beq\label{2.11}
m(r_j,R_j)\le\frac{1}{\pi}\log\frac{R_j}{r_j}+c_2.
\eeq
Furthermore, we claim that
\beq\label{2.12}
m(r_j,R_j)\ge\frac{4}{3\pi}\log\frac{R_j}{\sqrt{2}r_j}.
\eeq
Indeed, according to the nonnegativity of the module of a family of curves,
 there is no loss of generality in assuming that $R_j>\sqrt{2}r_j$.
Since
$$
|\ga(r)|\le\frac{3}{4}\pi r,\quad  \sqrt{2}r_j<r<R_j,
$$
we have
\begin{eqnarray*}
m(r_j,R_j)&\ge& m(\sqrt{2}r_j,R_j)
\ge m(\{\ga(r):\sqrt{2}r_j<r<R_j\})\\
&\ge&
\int_{\sqrt{2}r_j}^{R_j}\frac{dr}{|\ga(r)|}
\ge\frac{4}{3}
\int_{\sqrt{2}r_j}^{R_j}\frac{dr}{\pi r}
= \frac{4}{3\pi}\log\frac{R_j}{\sqrt{2}r_j}\, .
\end{eqnarray*}
Making use of (\ref{2.11}) and (\ref{2.12}), we obtain
$$
\frac{R_j}{r_j}\le 4\exp(3\pi c_2),
$$
which implies (\ref{2.10}).

\hfill$\Box$

\begin{lem}\label{lem2.3}
For $z\in\He$ with $|z-x_0|< 2c$,
\beq\label{2.13}
|F(z)-\eta_0|\le c_4 |z-x_0|,\quad
c_4:=\frac{c_1+1}{c}\exp(\pi(c_2+4+2\pi c_3^2)).
\eeq
\end{lem}
{\bf Proof}. 
First, we claim that for $w:=F(z)$,
\beq\label{2.n1}
|w-\eta_0|\le |w_0-\eta_0|.
\eeq
To prove (\ref{2.n1}), we assume that $|w-\eta_0|> |w_0-\eta_0|$.
Then, for the module of $\Ga_2':=\Ga(\eta_0,w;w_0,\infty;G)$ we have
\beq\label{2.n2}
m(\Ga_2')\le 4.
\eeq
Indeed, since for $\ga\in\Ga_2'$,
$$
\ga\cap\{\xi\in G:|\xi-\eta_0|=|w_0-\eta_0\}\neq\emptyset,
$$
the metric
$$
\rho_1(\xi):=\left\{\begin{array}{ll} \pi^{-1}&\mb{ if }
0\le\Re(\xi)\le\pi, |\Im(\xi-w_0)|\le 2\pi ,\\[2ex] 0&\mb{
elsewhere in $\C$}\end{array}\right.
$$
is admissible for $\Ga_2'$. Therefore,
$$
m(\Ga_2')\le\int\rho_1(\xi)^2dm(\xi)=4.
$$
This proves (\ref{2.n2}).

Moreover, 
according to our assumption
$|z-x_0|< e^{-4\pi}|z_0-x_0|$ and
the left-hand side of (\ref{2.8}), for the module of
$
\Ga_2:=F^{-1}(\Ga_2')=\Ga(x_0,z;z_0,\infty;\He)
$ we obtain
\beq\label{2.14}
m(\Ga_2)\ge\frac{1}{\pi}\log\left|\frac{z_0-x_0}{z-x_0}\right|
>4,
\eeq
which contradicts (\ref{2.n2}).
Hence, (\ref{2.n1}) holds.

Our next objective is to estimate from above the module of $\Ga_2'$ 
under the assumption (\ref{2.n1}).
Denote by $r(w)$ the supremum of values of $r>0$ such that $\ga(r)$ separates 
$\eta_0$ and $w$ in $G$. By Lemma \ref{lem2.2}
\beq\label{2.15}
|w-\eta_0|\le r(w)\le c_3|w-\eta_0|.
\eeq
Let $\rho_2$ be the extremal metric for the family $\Ga_3:=\Ga\left(
r(w),|\eta_0-w_0|\right)$,
i.e.,
$$
\int_\ga \rho_2(\xi)|d\xi|\ge 1,\quad \ga\in \Ga_3,
$$
$$
m\left(
r(w),|\eta_0-w_0|\right)=\int\rho_2(\xi)^2dm(\xi).
$$
Since for $\ga\in\Ga_2'$ with $\ga\cap\ga(r(w))\neq\emptyset$ we have
$|\ga|\ge r(w)/c_3$, for such $\ga$ and the metric
$$
\rho_{3}(\xi):=\left\{\begin{array}{ll} \frac{c_3}{r(w)}&\mb{ if }
\xi\in G, |\xi-\eta_0|\le 2r(w),\\[2ex] 0&\mb{
elsewhere in $\C$}\end{array}\right.$$
we obtain
$$
\int_\ga \rho_3(\xi)|d\xi|\ge 1.
$$
From what we already proved, we conclude that
the
metric
$$
\rho(\xi):=\max_{k=1,2,3}\rho_k(\xi),\quad \xi\in \C,
$$  
 is admissible for $\Ga_2'$ and, by (\ref{2.9}) and (\ref{2.15}),
\begin{eqnarray*}
m(\Ga_2')&\le&\sum_{k=1}^3\int\rho_k(\xi)^2dm(\xi)
\le4+\frac{1}{\pi}\log\frac{|w_0-\eta_0|}{r(w)} +c_2+
2\pi c_3^2\\
&\le& 
\frac{1}{\pi}\log\frac{|w_0-\eta_0|}{|w-\eta_0|} +c_5,\quad
c_5:=c_2+4+2\pi c_3^2.
\end{eqnarray*}
Comparing the above inequality with (\ref{2.14}) and 
 using (\ref{2.17}),
we obtain 
$$
|w-\eta_0|\le\frac{|w_0-\eta_0|}{|z_0-x_0|}e^{c_5\pi}|z-x_0|<
\frac{c_1+1}{c}e^{c_5\pi}|z-x_0|,
$$
which implies (\ref{2.13}).

\hfill$\Box$

\absatz{Proof of Theorem \ref{th1}}

Neither the hypothesis nor the conclusion is affected if
 we assume that $E\subset[-1,1],
\pm1\in E$ and $-1<x_0<1$. Next, we construct a sequence 
$E_j,j\in\N$ of compact sets each of which consists of 
a finite number of intervals such that $x_0\in$ Int$(E_j)$ for
each $j$ and
$$
E_{j+1}\subset E_j\subset[-1,1],\quad \bigcap_{j=1}^\infty E_j=E.
$$
By \cite[p. 108, Theorem 4.4.6]{ran}
\beq\label{3.1}
\lim_{j\to\infty} g_{\CC\setminus E_j}(z)=g_{\CC\setminus E}(z),
\quad z\in\C\setminus E.
\eeq
Since 
$
\En(|x-x_0|^\al,E)\le \En(|x-x_0|^\al,E_j),
$
by  assumption (\ref{1.1}) and the Vasiliev-Totik theorem (\ref{2.3}),
it follows that for some $\al>0$,
$$
h_{E_j}(x_0)^{-\al}\sigma_\al=
\lim_{n\to\infty}n^\al\En(|x-x_0|^\al,E_j)\ge
\limsup_{n\to\infty}n^\al\En(|x-x_0|^\al,E)=:\be>0.
$$
Hence,
$$
h_{E_j}(x_0)\le \left(\frac{\sigma_\al}{\be}\right)^{1/\al},\quad
\kap(E_j)\ge\kap(E),
$$
and (\ref{2.2}) holds for each $E_j$ instead of $E$ with a constant $c>2$
independent of $j$.
Referring to Lemma \ref{lem2.3}, we find that for $z\in \C\setminus E_j$ with
$|z-x_0|\le 2c$,
$$
g_{\CC\setminus E_j}(z)\le c_4|z-x_0|,
$$
 where a constant $c_4$ also does not depend on $j$. A passage to the limit
  as $j\to\infty$
and (\ref{3.1}) yield
\beq\label{3.2}
g_{\CC\setminus E}(z)\le c_4|z-x_0|.
\eeq
Next, consider the case  where $|z-x_0|\ge 2c>4$. By \cite[p. 53]{saftot},
in this case 
$$
g_{\CC\setminus E}(z)=\int\log|z-x|d\mu_E(x)-\log\kap(E)\le\log(2c|z-x_0|)
$$
which, together with (\ref{3.2}), implies (\ref{1.2})

This concludes the proof of Theorem \ref{th1}.

 \hfill$\Box$

\absatz{Acknowledgements}

 The
author is grateful to  M. Nesterenko
  for his helpful comments.

V. V. Andrievskii

 Department of Mathematical Sciences

 Kent State University

 Kent, OH 44242

 USA

e-mail: andriyev@math.kent.edu


\begin{thebibliography}{99}

\bibitem{ahl}
 Ahlfors L. V. (1973)

Conformal invariants, McGraw-Hill, New York.
  

\bibitem{and01}
   Andrievskii V. V. (2001)
	
	A Remez-type inequality in terms of capacity,
	Complex Variables 45, 35--46.


\bibitem{and08}
 Andrievskii V. V. (2008)

 Positive harmonic functions on Denjoy domains in
the complex plane,  J. d'Analyse Math.  104,  83--124.

\bibitem{and09}
Andrievskii V. V.  (2009)

Polynomial approximation of piecewise analytic
functions on a compact subset of the real line, J. Approx. Theory
161, 634--644.



\bibitem{andbla}  Andrievskii V. V. and  H. -P. Blatt (2002)

Discrepancy of signed measures and
      polynomial approximation,
       Springer-Verlag, Berlin/New York.


 \bibitem{ber1} Bernstein S. (1914)

 Sur la meilleure approximation de $|x|$ par
 des polynomes des degr\'es donn\'es, Acta Math. 37,  1--57.

 \bibitem{ber2} 
 Bernstein S.  (1938)

On the best approximation of $|x|^p$
 by means of polynomials of extremely high degree, Izv. Akad. Nauk
 SSSR Ser. Mat. 2,  160--180 (in Russian).

 \bibitem{ber3} 
 Bernstein S.  (1938)

On the best approximation of
 $|x-c|^p$, Dokl. Akad. Nauk SSSR 18,  379--384 (in Russian).


  \bibitem{cartot}
 Carleson L. and  V. Totik  (2004)

 H\"older continuity of Green's
 functions,  Acta Sci. Math. (Szeged)  70,     557--608.

 \bibitem{cargar}
 Carroll T. and  S. Gardiner  (2008)

  Lipschitz continuity of the
 Green function in Denjoy domains, Ark. Mat. 46,
271--283.


\bibitem{pom}  
  Pommerenke Ch. (1992)
	
	Boundary behavior of conformal maps, Vandenhoeck and Ruprecht,
                G\"ottingen.
            


\bibitem{ran}
  Ransford T. (1995)
		
		Potential theory in the complex plane,
   Cambridge University Press, Cambridge.


\bibitem{saftot}
 Saff E. B. and   V. Totik (1997)

    Logarithmic potentials with external
fields, Springer-Verlag,  New York/Berlin.



\bibitem{tot06}
 Totik V. (2006)

 Metric properties of harmonic measures,
Mem. Amer. Math. Soc. 184(867).


\bibitem{vas}
 Vasiliev R. K. (1998)

 Chebyshev polynomials and approximation theory on
 compact subsets of the real axis, Saratov University Publishing
 House, Saratov.

\bibitem{wid}
 Widom H. (1969)

 Extremal polynomials associated with a
system of curves in the complex plane,
 Adv. in Math.  3, 127--232.

\end{thebibliography}
\end{document}